\documentclass[10pt]{amsart}
\usepackage{mdwlist}
\usepackage{psfrag}

\usepackage{amsmath,amssymb,amsthm,latexsym}
\usepackage{verbatim}
\usepackage[dvips]{graphicx}

\newtheorem{thm}{Theorem}[section]
\newtheorem{lem}[thm]{Lemma}
\newtheorem{cor}[thm]{Corollary}
\newtheorem{conj}[thm]{Conjecture}
\newtheorem{prop}[thm]{Proposition}
\newtheorem{DEF}[thm]{Definition}
\theoremstyle{remark}
\newtheorem{rem}[thm]{Remark}

\newcommand{\thom}{\text{\tt Hom}\,}
\newcommand{\Homdm}{\texttt{Hom}\,(G,K_n)}
\newcommand{\Hom}{$\texttt{Hom}\,(G,K_n)$ }
\newcommand{\Homk}{$\texttt{Hom}\,(G,K_n)$}

\newcommand{\Hnulak}{$\texttt{Hom}_{\,0}(G,H)$}
\newcommand{\HomH}{$\texttt{Hom}\,(G,H)$ }

\newcommand{\vgap}{\text{\rm vgap}}
\newcommand{\conn}{\text{\rm conn}\,}
\newcommand{\nin}{\noindent}
\newcommand{\maxval}{\text{\rm maxval}}

\newcommand{\pet}{\vspace{5pt}}

\newcommand{\ub}{\textnormal{\bf u}}
\newcommand{\vb}{\textnormal{\bf v}}
\newcommand{\wb}{\textnormal{\bf w}}
\newcommand{\qb}{\textnormal{\bf q}}

\newcommand{\bdm}{\begin{displaymath}}
\newcommand{\edm}{\end{displaymath}}
\newcommand{\beq}{\begin{equation}}
\newcommand{\eeq}{\end{equation}}
\newcommand{\tn}{\textnormal}
\newcommand{\bn}{\partial}

\newcommand{\ahat}{A_1,\dots,\hat{A_j},\dots,A_p}

\title{Higher connectivity of graph coloring complexes}

\author{Sonja Lj. \v{C}uki\'{c}}
\address{Department of Mathematics, Royal Institute of Technology, Stockholm, Sweden}
\email{cukic@math.kth.se}
\author{Dmitry N. Kozlov}
\address{Institute of Theoretical Computer Science, ETH-Z\"urich, CH-8092 Z\"urich,
Switzerland, and Department of Mathematics, Royal Institute of Technology, Stockholm, Sweden (on leave)\newline \indent{\it E-mail addresses: }{\tt dkozlov@inf.ethz.ch, kozlov@math.kth.se} }
\thanks {Research is partially supported by the Swiss National Science 
Foundation, and by the Swedish Research Council}

\date{September 22, 2004}

\begin{document}

\begin{abstract}
The main result of this paper is a~proof of the following conjecture of Babson \& Kozlov:

\begin{quote}
{\bf Theorem.} {\em Let $G$ be a~graph of maximal valency~$d$, then the complex \emph{\Hom} is at least $(n-d-2)$-connected.}
\end{quote}

\nin Here $\thom(-,-)$ denotes the polyhedral complex introduced by Lov\'asz to study the topological lower bounds for chromatic numbers of graphs.

We will also prove, as a corollary to the main theorem, that the complex $\thom(C_{2r+1},K_n)$ is $(n-4)$-connected, for $n\geq 3$.
\end{abstract}

\maketitle

\section{Introduction.}

Given a~graph $G$ and a~positive integer $n$, there is a~construction
of a~topological space, which can be thought of as a~space of all
$n$-colorings of~$G$. This space, denoted $\thom(G,K_n)$, is
a~polyhedral complex, whose set of vertices coincides with the set of
all allowed vertex colorings of $G$ that use at most $n$ colors.

If $n<\chi(G)$, then this space is empty, otherwise, it possesses
a~meaningful topology. Intuitively, the cells of $\thom(G,K_n)$ of
positive dimension encode "homotopies of colorings", i.e., in some
sense, continuous procedures of changing one allowed coloring into
another one. The general definition is given in
subsection~\ref{ss2.2}.

$\thom(-,-)$-complexes were introduced by Lov\'asz, \cite{Lo1}, to study
topological lower bounds for chromatic numbers of graphs. The special
case $\thom(K_2,G)$ is the motivating example, since it turns out to
be homotopy equivalent to the~previously well-studied neighborhood
complex ${\tt N}(G)$. Recall that neighborhood complexes are the ones
which were so spectacularly used by Lov\'asz to resolve the Kneser
Conjecture, see \cite{Knes,Lo}.

Meanwhile, the other $\thom(-,-)$-complexes are not as
well-understood. The family of complexes $\thom(C_{2r+1},G)$
constitute one instance that has been studied lately, in connection with the proof by Babson \& Kozlov of the Lov\'asz Conjecture, \cite{BK03a,BK03b,BK03c}.

It has also been recently proven by the authors, see \cite{CK1}, that
for two arbitrary cycles $C_m$ and $C_n$, the connected components of
the complex $\thom(C_m,C_n)$ are either points, or are homotopy
equivalent to circles.

The intuition tells us that, if the maximal valency of the graph $G$
is small and the number $n$ is relatively large, then there should be
a~lot of freedom in completing the $n$-colorings of $G$
locally. Expressed topologically, we may hope that the complexes
$\thom(G,K_n)$ will be highly connected.

To say that $\thom(G,K_n)$ is $(-1)$-connected is the same as to say
that it is nonempty. It is well-known that a~sufficient condition is
provided by requiring that the maximal valency of $G$ is at most
$n-1$ (the most primitive coloring procedure works). Next, it was
shown in \cite[Proposition 2.4]{BK03b} that, if the maximal valency of
$G$ is at most $n-2$, then $\thom(G,K_n)$ is connected
($=0$-connected).

Motivated by these special cases and by some further computational
evidence, Babson \& Kozlov made the conjecture which initiated our
present study.

\begin{conj} \label{conj:main}
\cite[Conjecture 2.5]{BK03b}.

\nin The following inequality is valid for an arbitrary graph $G$:
\begin{equation} \label{eq:main}
 \conn\thom(G,K_n)\geq n-\maxval(G)-2.
\end{equation}
\end{conj}
\nin In other words: for an~arbitrary graph $G$, if $\maxval(G)< n-k$, then $\thom(G,K_n)$ is $(k-1)$-connected. 

The main purpose of this paper is to present a~proof of this conjecture. 
First, we analyze the situation for $n\geq\maxval(G)+3$. In this case, we find an~explicit algorithm for deforming an~arbitrary loop inside $\thom(G,K_n)$  to a~point, thereby verifying the triviality of the fundamental group.

Not surprisingly, it is virtually impossible to generalize these explicit homotopies to higher dimensions. Instead, we choose to use Hurewicz theorem, and calculate the nullity of the appropriate homology groups instead. Again, we give an~explicit algorithmic procedure for reducing an~arbitrary cycle to zero, by means of adding to it appropriately chosen boundaries.

\pet

{\bf Acknowledgments.} The first author would like to thank Rade
\v{Z}ivaljevi\'{c} for invaluable discussions and constant encouragement.

\section{Basic notations and definitions.}

%\deset

\subsection{Notations and terminology.} $\,$

\pet

\nin For an arbitrary natural number $n$ we introduce a shorthand
notation: $[n]=\{1,2,\dots,n\}$.

For any graph $G$, we denote the set of its vertices by $V(G)$, and
the set of its edges by $E(G)$, where $E(G)\subseteq V(G) \times
V(G)$. In this paper we will consider only undirected graphs, so
$(x,y)\in E(G)$ implies that $(y,x)\in E(G)$. Also, our graphs are
finite and \emph{do not} contain loops.

Let $G$ be a graph and $S\subseteq V(G)$. Then
we denote by $G[S]$ the graph with $V(G[S])=S$ and $E(G[S])=(S\times
S)\cap E(G)$. The graph $G[V(G)\setminus S]$ we will denote by $G-S$.

For a graph $G$ and $v\in V(G)$, let $\texttt{N}(v)=\{w\in V(G)\mid
(v,w)\in E(G)\}$, $\texttt{N}(v)$ is the set of neighbors
of~$v$. Denote the \emph{valency} of $v$ by $d(v)$, clearly
$d(v)=\vert \texttt{N}(v) \vert$. Also denote \emph{the maximal
valency} of $G$, $\max\{d(v)\mid v\in V(G)\}$, by $\maxval(G)$. 

For a graph $G$ and a natural number $n$, let $\vgap(G,n)$ denote what
we call the {\em $n$-color-valency gap}, namely, set
$\vgap(G,n):=n-\maxval(G)-1$. With these notations, the
equation~\eqref{eq:main} can be rewritten as: \[ \conn\thom(G,K_n)\geq
\vgap(G,n)-1.
\]
In other words, the first possibly nontrivial homotopy group of $\thom(G,K_n)$ is indexed by the $n$-color-valency gap, $\vgap(G,n)$.

Denote with $K_n$, and $C_n$ the complete graph (no loops), resp.\ the cycle with $n$ vertices, i.e., $V(K_n)=V(C_n)=[n]$, and
$E(K_n)=\{(x,y)\mid x,y\in[n],\ x\ne y\}$, $E(C_n)=\{(x,x+_n
1),(x+_n 1,x)\mid x\in[n]\}$. The maximal valency
of $C_n$ is 2, for $n \ge 3$, while the maximal valency of $K_n$ is
$n-1$.

For any two graphs $G$ and $H$, let $G\coprod H$ denote the disjoint
union of these graphs.

%\deset
\subsection{\texttt{Hom} complexes: definition, examples and basic 
properties} $\,$ \label{ss2.2}

\vspace{5pt}

\nin A standard generalization of graph colorings is provided by the
following definition.

\begin{DEF} 
For two graphs $G$ and $H$, a \textbf{graph homomorphism} from $G$ to
$H$ is a map $\phi: V (G)\to V(H)$ such that if $(x, y)\in E(G)$,
then $(\phi(x),\phi(y))\in E(H)$. 
\end{DEF}

We denote the set of all homomorphisms from $G$ to $H$ by
\textnormal{\Hnulak}.

\begin{DEF} \cite{Lo1},\cite[Definition 1.2]{BK03a}

\nin \textnormal{\HomH} is a polyhedral complex whose cells are indexed by all functions $\eta : V (G) \to 2^{V (H)}\setminus\{\varnothing\}$, such that if $(x, y)\in E(G)$, then for all $\tilde{x}\in \eta(x)$ and $\tilde{y}\in \eta(y)$, $(\tilde{x},\tilde{y})\in E(H)$. 

The closure of a cell $\eta$ consists of all cells indexed by $\tilde{\eta}: V(G)\to 2^{V(H)}\setminus\{ \varnothing\}$ which satisfy the condition that $\tilde{\eta}(v)\subseteq {\eta}(v)$, for all $v\in V(G)$.
\end{DEF}

\nin {\em Note.} We follow \cite{BK03a} in our notations.

The set of vertices of \HomH is \Hnulak. We note that cells of \HomH are direct products of
simplices, and that the dimension of a cell $\eta$ is equal to $\sum_{v \in
V(G)}\vert \eta (x) \vert -\vert V(G)\vert $.

One can describe a labeling of the cells in $\thom(G,K_n)$ rather
directly: they are indexed by all $p$-tuples of nonempty subsets of $[n]$,
$(A_1,\dots,A_p)$, where $p= |V(G)|$, and such that, if $(i,j)$ is an edge in $G$, then $A_i\cap A_j=\varnothing$. Also, for a cell $\eta\in \thom(G,K_n)$, we will denote with $A_j^\eta$ its $j$-th coordinate set (and sometimes we will refer to it as {\it $j$-th color list} of $\eta$). With these notations, $\eta=(A_1^\eta,\dots,A_p^\eta)$.

For any three graphs $G$, $H$ and $K$, the following is
true:
\begin{equation}
\label{disjunija}
\texttt{Hom}\,(G\coprod H,K)=\texttt{Hom}\,(G,K)\times
\texttt{Hom}\,(H,K).
\end{equation}

%\deset
%{\bf Examples} of \Hom complexes:\\
%The most of these examples were observed in \cite{BK03b}.
%\begin{itemize}
%\item $\texttt{Hom}\,(G,K_1)$ is empty set if $E(G)\neq \varnothing$; otherwise it is a point.
%\item $\texttt{Hom}\,(G,K_2)=\varnothing$ if $G$ is not bipartite; if $G$ is bipartite with $k$ connected components then $\texttt{Hom}\,(G,K_2)$ consists of $2^k$ points.
\begin{figure}[ht] 
\begin{center} 
\psfrag{Y}{$K_{1,3}$}
\psfrag{H}{$W_4$}
\psfrag{HomY}{$\texttt{Hom}\,(K_{1,3},K_3)$}
\psfrag{HomH}{$\texttt{Hom}\,(W_4,K_4)$}
%\psfrag{Hom(c5,k3)}{$\texttt{Hom}\,(C_5,K_3)$}
\includegraphics[scale=1]{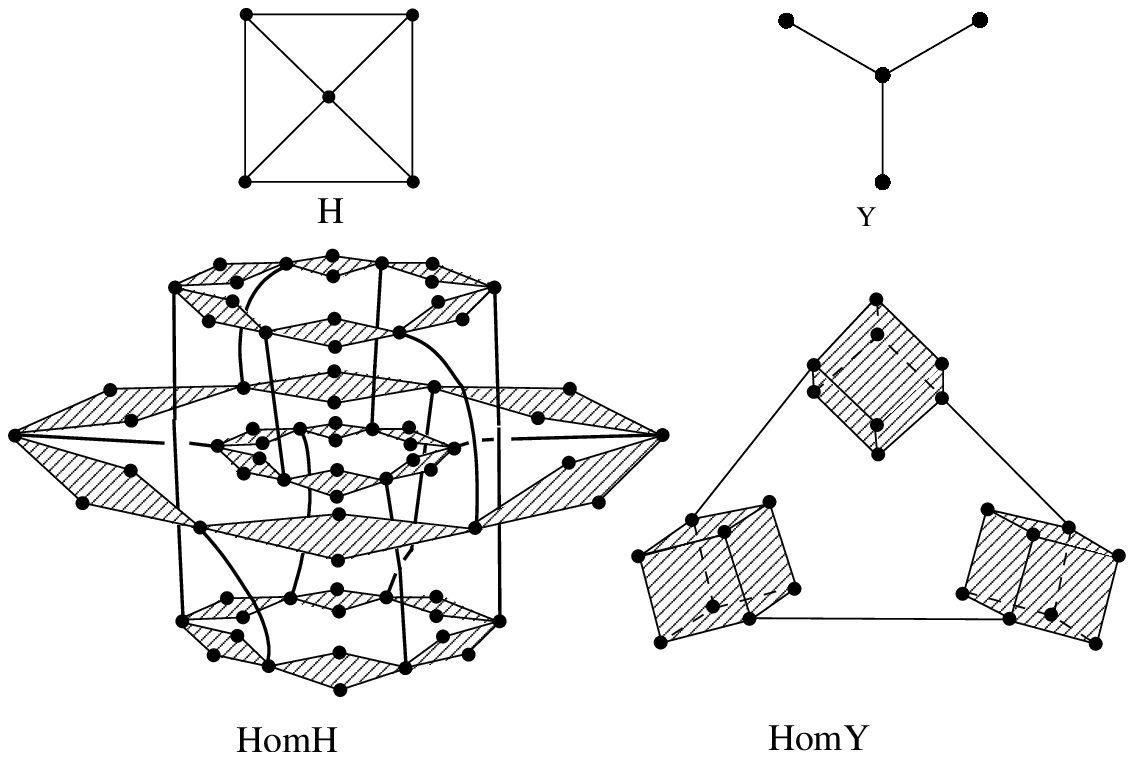}
\caption{$\,$} 
\label{slika}
\end{center} 
\end{figure}

Two examples of complexes $\thom(G,K_n)$ are shown on {\sc Figure}
\ref{slika}. Note that $\thom (K_{1,3},K_3)$ is a $3$-dimensional complex, with 3 solid cubes and 3 intervals. Further examples can be found in~\cite{BK03b}.

%\item $\texttt{Hom}\,(C_6,K_3)$ consists of six isolated points, 18 squares and 6 solid cubes.
%\item $\texttt{Hom}\,(C_7,K_3)$ is homeomorphic to a disjoint union of two M\"{o}bius bands.
%\item $\texttt{Hom}\,(K_2,K_4)$ is the boundary of a cuboctahedron (convex quasi-regular polyhedra with 14 faces: 8 triangles and 6 squares).
%\item $\texttt{Hom}\,(K_n,K_n)$ is a disjoint union of $n!$ points.
%\end{itemize}

\section{The fundamental group of {\textnormal \Hom}.}

%\deset

\nin We are now ready to study the fundamental
group of the graph coloring complexes. 

Let $\lambda(G)$, resp. $p(G)$, denote the cardinality of a~maximal independent set
in $G$, resp. the number of vertices of $G$. We shall write $\lambda$ and $p$ whenever it is clear which $G$
is meant. 
For an arbitrary graph $G$, label the vertices with $x_1,\dots,x_p$, such that  $\{x_1\dots,x_\lambda\}$ is a maximal independent set.

\begin{lem}\label{dokaz} 
Let $G$ be an arbitrary graph, and $n\geq 2$. Assume that
$\vgap(G,n)\geq 2$. For an~arbitrary $1<i\leq n$, any closed edge-path
in the $1$-skeleton of \emph{\Homk}, such that the first $\lambda$
coordinates of every vertex in this path are elements of the set
$\{i-1,i,\dots,n\}$, can be deformed by subsequent homotopies so that
these $\lambda$ coordinates will be elements of the set
$\{i,i+1,\dots,n\}$.

\begin{proof}
Each vertex of the complex \Hom can be described by $p$-tuple
$(a_1,a_2,\dots,a_p)$, where $a_d\in[n]$, for $d\in[p]$.

First of all, in any closed edge-path $u_1,u_2,\dots,u_s=u_1$ we may assume that successive vertices are distinct. Otherwise, if $u_d=u_{d+1}$ for some $d$, we could delete $u_d$ and obtain a homotopic path. 
Hence, from now on, we shall implicitly  replace every path with a homotopic path where no two successive  vertices are equal.

We will give an algorithm which performs  the deformation, the existence of which is claimed in our lemma. Formal description of steps, together with proofs of their correctness, will be given after description of the algorithm.

Let $\ub^{(\lambda+1)}$ be  a~closed path satisfying conditions of Lemma~\ref{dokaz}, and let $u_1^{(\lambda+1)},u_2^{(\lambda+1)},\dots,u_m^{(\lambda+1)}=u_1^{(\lambda+1)}$ be the vertices of this path. Note that each $u_k^{(\lambda+1)}$, for $k\in[m-1]$, is described by a $p$-vector and  $u_k^{(\lambda+1)}$ differs from  $u_{k+1}^{(\lambda+1)}$ in exactly one component.
\begin{figure}[ht] 
\begin{center}
\psfrag{1}{ $1$} 
\psfrag{t}{ $\vdots$} 
\psfrag{i-2}{ $i-2$} 
\psfrag{i-1}{ $i-1$} 
\psfrag{i}{ $i$} 
\psfrag{n}{ $n$}
\psfrag{x1}{ $x_1$} 
\psfrag{ht}{ $\cdots$} 
\psfrag{l}{ $x_\lambda$} 
\psfrag{l+1}{ $x_{\lambda+1}$} 
\psfrag{p}{ $x_p$} 
\psfrag{colors}{\it{Colors}} 
\psfrag{vertices}{\it{Vertices of $G$}} 
\psfrag{12}{\sc{Steps 1 and 2}} 
\psfrag{34}{\sc{Steps 3 and 4}} 
%\psfrag{i-2}{ $i-2$} 
 \includegraphics[scale=0.9]{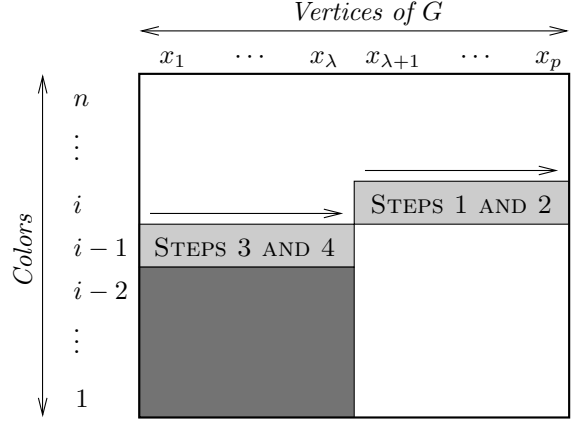}
\caption{Scheme of the proof of Lemma \ref{dokaz}: shaded areas indicate ``forbidden colors''; before we ``get rid'' of color $i-1$ on first $\lambda$ vertices of $G$, we eliminate color $i$ from the vertices $x_{\lambda+1},\dots,x_p$.} 
\label{scheme}
\end{center} 
\end{figure}

The input for our algorithm is the path $\ub^{(\lambda+1)}$, and the algorithm consists of two parts. In the first part, we repeat both of the following steps for all $j=\lambda+1,\dots,p$, in increasing order.
\begin{description}
\item[Step 1] Inserting two new vertices between all those neighboring pairs $u_{k}^{(j)}$ and $u_{k+1}^{(j)}$ from the path $\ub^{(j)}$ which both have $i$ on $j$-th position; result of step 1 is a path $\vb^{(j)}$.
\item[Step 2] Deleting those vertices from the path  $\vb^{(j)}$ which use color $i$ on the graph vertex $x_j$; result of this step is the path  $\ub^{(j+1)}$. 
\end{description} 
After the first phase, we obtain a path $\ub^{(p+1)}$ with the property that colors used on vertices $x_{\lambda+1},\dots,x_p$ are elements of the set $[n]\setminus\{i\}$.

Now, input for the next stage is the path $\wb^{(1)}=\ub^{(p+1)}$, and we repeat both steps 3 and 4 for all $j=1,\dots,\lambda$, again in the increasing order. 
\begin{description}
\item[Step 3] Inserting two new vertices between those neighboring vertices of the path $\wb^{(j)}$,  $w_{k}^{(j)}$ and $w_{k+1}^{(j)}$, which both have $i-1$ on $j$-th position; result of step 3 is the path $\qb^{(j)}$.
\item[Step 4] Deleting those vertices from the path  $\qb^{(j)}$ which use color $i-1$ on $x_j$; result of this step is the path  $\wb^{(j+1)}$. 
\end{description} 
Output of this algorithm is a path $\wb^{(\lambda+1)}$, and for each vertex of this path, colors used on $x_1,\dots,x_\lambda$ are from the set $\{i,i+1,\dots,n\}$.\\

Now we give detailed description of all the steps.

\nin {\bf Step 1.} If $i$ occurs in $j$-th position for both $u_{k}^{(j)}$ and $u_{k+1}^{(j)}$ we will ``separate'' $u_{k}^{(j)}$ and $u_{k+1}^{(j)}$  by adding new vertices to the given path. We have
%It is obvious that we can assume that $u_{k}^{(j-1)}\ne u_{k+1}^{(j-1)}$ (otherwise we could delete one of them and obtain homotopic path).
\bdm \begin{split}
u_{k}^{(j)}&=(\dots,x,\dots,\overbrace{\quad i\quad}^{j\textrm{\tiny -th pos}},\dots)\\
u_{k+1}^{(j)}&=(\dots,y,\dots,\,\quad {i}\quad \, ,\dots)
\end{split}\edm
%\end{center}
where $x\ne y$. We have indicated only $i$ and components which are different in these two vertices. Let 
\[ 
A=\{\mu\textnormal{-th coordinates of } u_{k}^{(j)}\} \cup\{\mu\textnormal{-th coordinates of } u_{k+1}^{(j)}\},
\]
where both sets are indexed with those $\mu$, for which $x_\mu \in \texttt{N}(x_j)$. Obviously, we have $|A\cup\{i\}|\leq\maxval(G)+2$ and, since $n\ge\maxval(G)+3$, there exists $z\in [n]\setminus(A\cup \{i\})\neq\emptyset$, such that
%\begin{center}
\bdm \begin{split}
U_{k}^{(j)}&=(\dots,x,\dots,\overbrace{\quad z \quad}^{j\textrm{\tiny -th pos}},\dots)\\
U_{k+1}^{(j)}&=(\dots,y,\dots,\quad {z}\quad \,,\dots)
\end{split} \edm
%\end{center}
are vertices of \Homk. In another words, we have a vertex of $G$ of bounded degree, and two colorings that differ in only one vertex, so we find an alternative color left for the vertex $x_j$.

 Deformation of the subpath  $u_{k}^{(j)}u_{k+1}^{(j)}$ into  $u_{k}^{(j)}U_{k}^{(j)}U_{k+1}^{(j)}u_{k+1}^{(j)}$ is a homotopy over the 2-cell $(\dots,\{x,y\},\dots,\{z,i\},\dots)$, see Figure \ref{add_vert}.

\begin{figure}[ht] 
\begin{center}
\psfrag{1}{ $u_{k}^{(j)}$} 
\psfrag{2}{ $u_{k+1}^{(j)}$} 
\psfrag{3}{ $u_{k}^{(j)}$} 
\psfrag{4}{ $u_{k+1}^{(j)}$} 
\psfrag{5}{ $U_{k+1}^{(j)}$} 
\psfrag{6}{ $U_{k}^{(j)}$} 
\includegraphics[scale=1.0]{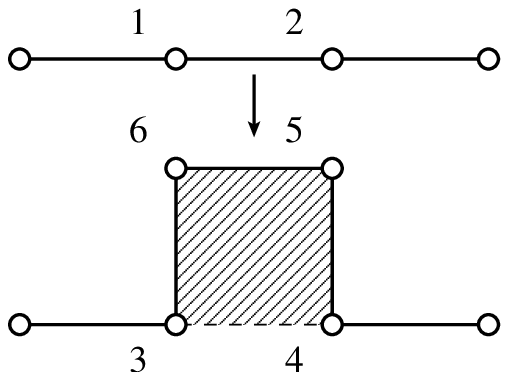}
\caption{$\,$} 
\label{add_vert}
\end{center} 
\end{figure}

\nin Let $\vb^{(j)}=v_1^{(j)},v_2^{(j)},\dots,v_1^{(j)}$ be the path obtained by the first step.

\pet

\nin {\bf Step 2.} In this step we will remove from the path $\vb^{(j)}$ those vertices which have $i$ in the $j$-th position. In this case, as a result of the first step, we have the following situation: 
%\begin{center}
\bdm \begin{split}
v_{k-1}^{(j)}&=(\dots,\overbrace{\quad x \quad}^{j\textrm{\tiny -th pos}},\dots)\\
v_{k}^{(j)}&=(\dots,\,\quad {i}\quad \, ,\dots)\\
v_{k+1}^{(j)}&=(\dots,\quad y \quad\,,\dots)
\end{split} \edm
%\end{center}
where $x,y\ne i$.

Deletion of $v_{k}^{(j)}$ is either a~homotopy over the 2-cell
$(\dots,\{x,y,i\},\dots)$, or, if $x=y$, a~homotopy over the
1-cell $(\dots,\{x,i\},\dots)$, see Figure \ref{rem_vert}.

\begin{figure}[ht] 
\begin{center}
\psfrag{1}{ $v_{k-1}^{(j)}$} 
\psfrag{2}{ $v_{k+1}^{(j)}$} 
\psfrag{3}{ $v_{k}^{(j)}$} 
\includegraphics[scale=1.0]{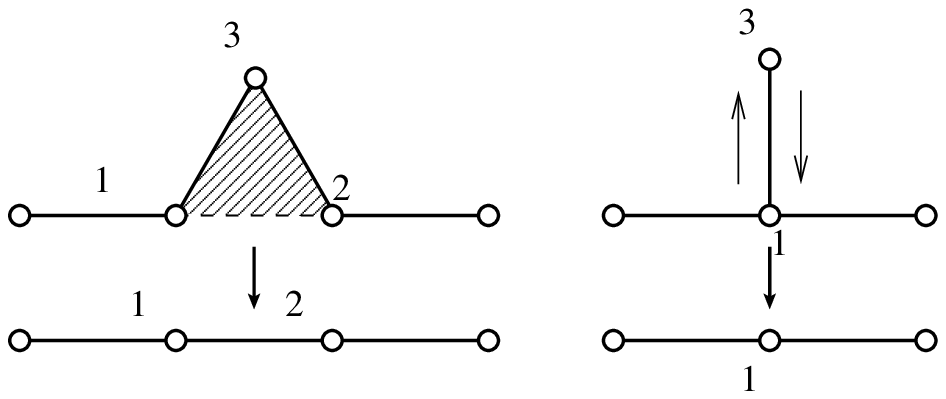}
\caption{$\,$} 
\label{rem_vert}
\end{center} 
\end{figure}

Let $\ub^{(j+1)}=u_1^{(j+1)},u_2^{(j+1)},\dots,u_1^{(j+1)}$ be the path obtained after
this step. It is clear that this path does not have any vertices
which have $i$ in the $j$-th position.

%\suspend{enumerate}
After the first stage of our algorithm, we obtain a path $\ub^{(p+1)}$. Set $\wb^{(1)}=\ub^{(p+1)}$, that is $\wb^{(1)}=w_1^{(1)},w_2^{(1)},\dots,w_1^{(1)}$, where $w_d^{(1)}=u_d^{(p+1)}$ for all $d\in [l-1]$ ($l$ is the length of the path $\ub^{(p+1)}$). In this path, $i$ is {\it not} in $j$-th position of any vertex, where $j\in\{\lambda+1,\dots,p\}$.

%Now, we will repeat two similar steps for all $j\in\{1,2,\dots,\lambda\}$ in the increasing order of~$j$.
%\resume{enumerate}

\nin {\bf Step 3.} If $i-1$ occurs in $j$-th position for both $w_{k}^{(j)}$ and $w_{k+1}^{(j)}$ we will, similar to the first step, ``separate'' them by adding new vertices to the given path. We start with
%\begin{center}
\bdm \begin{split}
w_{k}^{(j)}&=(\dots,\overbrace{i-1}^{j\textrm{\tiny -th pos}},\dots,x,\dots)\\
w_{k+1}^{(j)}&=(\dots,\ {i-1}\ ,\dots,y,\dots)
\end{split} \edm
%\end{center}
where $x\ne y$. Since vertices labeled $x_1,\dots,x_\lambda$ form an independent set, we know that ${\tt N}(x_j)\subseteq\{x_{\lambda+1},\dots,x_p\}$. Hence, the color $i$ does not occur among the neighbors of $x_j$ in $w_k^{(j)}$, or in $w_{k+1}^{(j)}$. Therefore,
%\begin{center}
\bdm \begin{split}
W_{k}^{(j)}&=(\dots,\overbrace{\quad i\quad}^{j\textrm{\tiny -th pos}},\dots,x,\dots)\\
W_{k+1}^{(j)}&=(\dots,\,\quad {i}\quad \,,\dots,y,\dots)
\end{split} \edm
%\end{center}
are legal $n$-colorings, and hence are also vertices of $\thom(G,K_n)$.
Then, the deformation of the subpath  $w_{k}^{(j)}w_{k+1}^{(j)}$ into  $w_{k}^{(j)}W_{k}^{(j)}W_{k+1}^{(j)}w_{k+1}^{(j)}$ is a~homotopy over the 2-cell $(\dots,\{i-1,i\},\dots,\{x,y\},\dots)$.

Let now $\qb^{(j)}=q_1^{(j)},q_2^{(j)},\dots,q_1^{(j)}$ be the path obtained by the third step.

\nin {\bf Step 4.} Similar to the Step 2, all the vertices from the path $\qb^{(j)}$ which have $i-1$ in the $j$-th position are deleted:
%\begin{center}
\bdm \begin{split}
q_{k-1}^{(j)}&=(\dots,\overbrace{\quad x \quad}^{j\textrm{\tiny -th pos}},\dots)\\
q_{k}^{(j)}&=(\dots,\; i-1 \;,\dots)\\
q_{k+1}^{(j)}&=(\dots,\quad y\quad\,,\dots)
\end{split} \edm 
%\end{center}
where $x,y\ne i-1$.

Deletion of $q_{k}^{(j)}$ is homotopy over the 2-cell $(\dots,\{x,y,i-1\},\dots)$ (or, if $x=y$, over 1-cell $(\dots,\{x,i-1\},\dots)$). 

Let $\wb^{(j+1)}=w_1^{(j+1)},w_2^{(j+1)},\dots,w_1^{(j+1)}$ be the obtained path. We can
see that after this step $j$-th coordinate of any vertex from this
path is in the set $\{i,\dots,n\}$.\\

Finally, we arrive at the path $\wb^{(\lambda+1)}$ which has been obtained by a~sequence of elementary homotopies from the original path and has the additional property that the first $\lambda$ coordinates of every vertex from this path are elements of the set $\{i,\dots,n\}$.
\end{proof}
\end{lem}

\nin {\em Note.} This proof was motivated by the ideas from \cite{TopG}.

\pet

 We would like to remark, that the Steps 1 and 2 can be combined to obtain one reduction step encoding the following transformation: we find the first vertex where $i$ occurs in $j$-th position of two vertices in a~row, then we glue in a~square, as in Step~1, and then we clip off a~triangle, as in Step~2, see Figure~\ref{rstep}.

\begin{figure}[ht] 
\begin{center}
  \begin{picture}(0,0)%
    \includegraphics{red.pstex}%
  \end{picture}%
  \input{red.pstex_t}%
  
\caption{$\,$} 
\label{rstep}
\end{center} 
\end{figure}

 The outcome of this procedure will be a~shortening of the undesired part (here meaning $i$ is in the $j$-th position) by~1. With this line of argument, one has two special cases to attend to. First, if the undesired vertices come only as singletons, then we clip them all off as in Step~2. Second, if all vertices are undesired, then gluing in an~arbitrary square, as in Step~1, reduces this to the case which we considered first.

The Lemma~\ref{dokaz} is the crucial step in the proof of the main result of this chapter.
 
\begin{thm}\label{fgrupa} 
Let $G$ be any graph. If $\,\vgap(G,n)\geq 2$, then \textnormal{\Hom} is simply connected.

%\begin{Note} It was proven in \cite{BK03b} that \Hom is connected if $n\ge d+2$. \end{Note}

\begin{proof}
We will prove that \Hom is simply connected using induction on the
maximal valency of $G$. We will use the same notations as in the
previous lemma.

Suppose that the maximal valency of a graph $G$ is $0$ and that
$n\ge3$. Then $G$ is disjoint union of $p$ points and, using
(\ref{disjunija}), we conclude:
\begin{center}
$\texttt{Hom}\,(G,K_n)=\texttt{Hom}\,(K_1,K_n)\times\texttt{Hom}\,(K_1,K_n)\times\dots\times\texttt{Hom}\,(K_1,K_n)$.\end{center}
$\texttt{Hom}\,(K_1,K_n)$ is a simplex, so \Hom is contractible, and hence $1$-connected.

Assume now that the maximal valency of $G$ is equal to $d\geq 1$ and
that $\vgap(G,n)\geq 2$. Let $\beta_1,\beta_2,\dots,\beta_1$ be any closed
edge-path in $1$-skeleton of \Homk. Since \Hom is a polyhedral complex,
it is clear that it is sufficient to consider only these paths. 

Using Lemma~\ref{dokaz} iteratively, we can homotopicaly transform this path to a~path $\alpha_1,\alpha_2,\dots,\alpha_1$ which has the property that the first $\lambda$ coordinates of any vertex from this path are equal to $n$, where $\{x_1,x_2,\dots,x_\lambda\}$ is a~maximal independent set
in $G$. Hence, the original path is homotopic to a path lying inside
the subcomplex $\texttt{Hom}\,(G-\{x_1,x_2,\dots,x_\lambda\},K_{n-1})$. That path is contractible to a point by induction hypothesis, since $\{x_1,\dots,x_\lambda\}$ is a~maximal independent set, and therefore the maximal valency of $G-\{x_1,x_2,\dots,x_\lambda\}$ is strictly less then~$d$.
\end{proof}
\end{thm}

We will now evaluate the fundamental group of $\texttt{Hom}\,(G,K_{n+1})$ for one case when the maximal valency of a graph $G$ equals $n-1$.
\begin{prop} \textnormal {$\pi_1(\texttt{Hom}\,(K_{n},K_{n+1}))$} is the free product of $\alpha_n$ copies of $\mathbb{Z}$, where $\alpha_n=n!\frac{n^2-n-2}{2}+1$.

\begin{proof}
Since $\texttt{Hom}\,(K_{n},K_{n+1})$ is a connected graph, we can choose a~spanning tree and contract it. Then we get a bouquet of $e-v+1$ circles, where $v$ is the number of vertices and $e$ number of edges of this graph. It is easy to see that $v=(n+1)!$ and, since $\texttt{Hom}\,(K_{n},K_{n+1})$ is $n$-regular, $e=\frac{nv}{2}$. This gives us the claim of the proposition. 
\end{proof} 
\end{prop} 

\section{Homology groups of \textnormal{\Homk} complexes and the main theorem.}

In the previous section we proved that, if the maximal valency of $G$ is $d$, then \Hom is simply connected (or $1$-connected), for all $n\geq d+3$. Here we will prove a more general statement using the notations and ideas from Lemma \ref{dokaz} and  Theorem \ref{fgrupa}.

First, we will introduce a new notation. Let $G$ be a graph with the set of vertices $V(G)=\{x_1,\dots,x_p\}$ and let $1\leq j \leq p$ be an integer. Then, for all $i\in [n]$, we define a~subcomplex $X_i(G,j)$ of \Hom in the following way:
%\bdm
	\[X_i(G,j)=\{\eta \in \Homdm\mid\eta(x_q)\subseteq \{i,i+1,\dots,n\},\tn{ for all }q\in[j] \}.
\]
%\edm
In other words, only colors $i,\dots,n$ have been used in the first $j$ vertices. For example, Lemma \ref{dokaz} is exactly pushing the loops from $X_{i-1}(G,\lambda)$ into the subcomplex $X_i(G,\lambda)$. If it is clear which graph $G$ is meant, we will use the notation $X_i(j)$ instead of  $X_i(G,j)$.

Let us again label the vertices of $G$ in the same way as we did in Lemma~\ref{dokaz}, that is so that the vertices labeled $x_1,\dots,x_\lambda$ form a~maximal independent set.

\begin{lem}\label{dokaz_hom} Let $G$ be a graph with the maximal valency equal to $d\geq 1$. If $C$ is a $t$-cycle in $X_{i-1}(\lambda)$, where $1\leq t \leq n-d-2$ and $i\in\{2,3,\dots,n\}$, then there exist a $t$-cycle $C'$ in $X_{i}(\lambda)$ such that $C$ and $C'$ represent the same element in \textnormal{$H_t(\Homdm,\mathbb{Z})$}.

\begin{proof}
Recall that a cell $\eta$ from \Hom can be described by the $p$-tuple $(A_1^\eta,A_2^\eta,\dots,A_p^\eta)$, where $A_j^\eta\subseteq [n]$ for all $j\in[p]$ ($p=\vert V(G) \vert$). Rephrasing the definition above, we get
\[
X_q(\lambda)=\{\eta \in \Homdm \mid A_1^\eta,\dots,A_\lambda^\eta \subseteq \{q,q+1,\dots,n\}\},
\]
 for $q\in[n]$. The orientation of the cells $\eta \in \Homdm$ can be chosen so that the boundary operator $\bn$ is given by
\[
\bn(\eta)=\bn (A_1^\eta,\dots,A_p^\eta)=\sum(-1)^{c(x)}(A_1^\eta,\dots,A_{s(x)}^\eta\setminus\{x\},\dots,A_p^\eta),
\] 
where $s(x)$ denotes the index, for which $x\in A_{s(x)}^\eta$, $c(x)=\vert A_1^\eta\vert+\dots+\vert A_{s(x)-1}^\eta\vert + \vert\{y\mid y\in A_{s(x)}^\eta\tn{ and } y<x\}\vert$, and the sum is taken over all $x\in A_1^\eta\cup\dots\cup A_p^\eta$, such that $\vert A_{s(x)}^\eta\vert\geq 2$.

Like in the proof of Lemma \ref{dokaz}, we give a description of an algorithm whose output is a cycle $C'$, whose existence is claimed in this lemma. General scheme of the proof is the same as for already mentioned Lemma \ref{dokaz}. Namely, we first get rid of color $i$ on vertices $x_{\lambda+1},\dots,x_p$, and afterwards we are removing color $i-1$ from color lists on all vertices from the chosen maximal independent set, see Figure \ref{scheme}.

Input for our algorithm is the $t$-cycle $C_{\lambda+1}=C$. For each $j=\lambda+1,\dots,p$, in the increasing order of~$j$, we repeat both of the next two inductive steps.
\begin{description}
\item[Step 1] This is the only step in our proof where we use the assumption about maximal valency of $G$. In this step we eliminate all cells $\eta$ from the cycle $C_j$ such that $A_j^\eta=\{i\}$. Result of this step is a cycle $C^2_j$ with the property that $j$-th color list of each cell from $C^2_j$ which contains $i$ has length $\geq 2$. 
\end{description}
\hspace{17pt} \begin{tabular}{p{11.66cm}}Now we have new, iterative, step. Namely, step 2 is repeated for cycles $C_j^l$, starting from $l=2$, until there exists a cell $\sigma\in C_j^l$ such that $i\in A_j^\sigma$. \end{tabular}
\begin{description}
\item[Step 2] Transforming the cycle $C_j^l$ to a homologous cycle $C_{j}^{l+1}$ with the property that, for all $\eta\in C_j^{l+1}$, if $i\in A_j^\eta$, then $\vert A_j^\eta\vert\geq l+1$.
 \end{description}
\hspace{17pt} \begin{tabular}{p{11.66cm}}Since length of each color list is certainly less or equal to $n$, repetition of step 2 will stop after finite number of iterations, and we will get a cycle $C_{j+1}$.\end{tabular} \\

\nin Input for the second stage of this algorithm is the cycle $C_1=C_{p+1}$. For each cell $\sigma$ of this cycle, $i\notin A_{\lambda+1}^\sigma \cup \dots A_p^\sigma$. Again, we repeat the following inductive steps for all $j=1,\dots,\lambda$, in the increasing order:
\begin{description} \item[Step 3]  We eliminate all cells $\eta$ from the cycle $C_j$ such that $A_j^\eta=\{i-1\}$. Result of this step is a cycle $C^2_j$  with the property that $j$-th color list of each cell from $C^2_j$ which contains $i-1$ has length $\geq 2$.
\end{description}
\hspace{17pt} \begin{tabular}{p{11.66cm}} Step 2 is repeated for cycles $C_j^l$, starting from $l=2$, until there exists a cell $\sigma\in C_j^l$ such that $i-1\in A_j^\sigma$. \end{tabular}
\begin{description}
\item[Step 4] Transforming the cycle $C_j^l$ to a cycle $C_{j}^{l+1}$, such that $C_j^l$ and $C_{j}^{l+1}$ represent the same homology element, and such that,  if $i-1\in A_j^\eta$, then $\vert A_j^\eta\vert\geq l+1$, for all cells $\eta\in C_j^{l+1}$.
 \end{description}
\hspace{17pt} \begin{tabular}{p{11.66cm}} By the same argument as in the first part, repetition of step 4 will stop after finite number of iterations, and we will get a cycle $C_{j+1}$. In the case when $j=\lambda$, we label the resulting cycle with $C'$.\end{tabular} \\

After describing the algorithm, we give a detailed description of steps:\\
\nin {\bf Step 1.} Let $\eta$ be a $t$-cell appearing in the chain $C_j$, such that $A_j^\eta=\{i\}$, and let $\eta$ have the coefficient $k\in \mathbb{Z}$ in $C_j$. The dimension of $\eta$ is equal to $t$, and we have:
\bdm 
\begin{split}
t =\sum_{r=1}^p (\vert A_r^\eta\vert-1) &=(\sum_{r\vert x_r\in \texttt{N}(x_j)}\vert A_r^\eta\vert-\vert \texttt{N}(x_j)\vert)+(\vert A_j^\eta\vert-1)+\\ &+(\sum_{\substack{r\vert x_r\notin \texttt{N}(x_j),\\ r\neq j}}\vert A_r^\eta\vert-(p-\vert \texttt{N}(x_j)\vert-1)) 
\geq \sum_{\substack{r\in [p] \\ x_r\in \texttt{N}(x_j)}}\vert A_r^\eta\vert-d,
\end{split}
\edm 
since the second term in the middle sum is~0, the third term in the middle sum is nonnegative, and $|\texttt{N}(x_j)|\leq d$. It follows that 
	\[\sum_{\substack{r\in [p] \\ x_r\in \texttt{N}(x_j)}}\vert A_r^\eta\vert\leq t+d\leq n-2;
\]
 as we recall that, by assumption of the lemma, $t\leq n-d-2$. Hence, there exists $I(\eta)\in [n]\setminus \{i\}$ such that $\eta'=(A_1^\eta,\dots,A_{j-1}^\eta,\{i,I(\eta)\},A_{j+1}^\eta,\dots,A_p^\eta)\in \Homdm$. Obviously, $\dim \eta'=t+1$. Set $C'_{j}=C_j-(-1)^{K(\eta')}k \ \bn \eta'$, where $(-1)^{K(\eta')}$ comes from the appropriate incidence number, i.e., in our notations
\bdm
K(\eta')=\vert A_1^\eta \vert + \dots+ \vert A_{j-1}^\eta \vert+\left\{ \begin{array}{ll} 1, & \tn{ if } i<I(\eta); \\ 0, & \tn{otherwise.} \end{array} \right.
\edm  
Clearly, $C'_{j}$ represents the same homology element as $C_j$. Furthermore, the number of cells $\sigma$ appearing in $C_j'$, such that $A_j^\sigma=\{i\}$, is strictly less then number of such cells appearing in $C_j$, since we have just eliminated the appearance of the cell~$\eta$. We repeat this procedure until we get cycle $C_j^2$ in which there are no cells  $\sigma$ such that $A_j^\sigma=\{i\}$.

\nin {\bf Step 2.} Suppose now that we have a cycle $C_j^l$, $l\geq 2$, which represents the same element in homology as $C_j$, and which has the additional property that $\vert A_j^\eta\vert\geq l$, for each cell $\eta$ appearing in $C_j^l$, such that $i\in A_j^\eta$. If no such cell exists, we are done with this case, and we set $C_{j+1}:=C_j^l$. Since $l\leq n$, we will always come to this case after a~finite number of steps. 

Assume now there exists a~cell $\eta$ such that $i\in A_j^\eta$, and let us construct $C_j^{l+1}$. Assume further that $C_j^l=\sum_{q\in Q} k_q\sigma_q$, where $Q\subset\mathbb{N}$, $k_q\in \mathbb{Z}$, $k_q\neq 0$, and $\sigma_q$'s are pairwise different
$t$-cells from \Homk. Then we can write $C_j^l=D_l+D_{>l}+D_0$, where
\begin{eqnarray*}
& &D_l=\sum_{q\in Q_l}k_q\sigma_q\tn{ where }Q_l=\{q\in Q\mid i \in
A_j^{\sigma_q} \tn{ and }\vert A_j^{\sigma_q} \vert=l\}, \\ &
&D_{>l}=\sum_{q\in Q_{>l}}k_q\sigma_q\tn{ where }Q_{>l}=\{q\in Q\mid i \in A_j^{\sigma_q} \tn{ and }\vert A_j^{\sigma_q} \vert>l\},\\ & &
D_0=\sum_{q\in Q_0}k_q\sigma_q\tn{ where }Q_0=\{q\in Q\mid i \notin
A_j^{\sigma_q}\}.
\end{eqnarray*}

If $D_l=0$, then we set $C_j^{l+1}=C_j^l$, so assume $D_l\neq 0$. 

Clearly, $\bn C_j^l=\bn D_l+\bn D_{>l}+\bn D_0$. Now let $\bn D_l=\beta+\gamma+\delta$, where $\beta$ is a subchain consisted of all cells $\sigma$ from the boundary such that $\vert A_j^\sigma \vert= l$, $\gamma$ is a chain of all cells from $\bn D_l$ such that $j$-th set has $l-1$ elements and contains $i$; finally $\delta$ consists of all cells where we do not have $i$ in $j$-th coordinate set. Since $0=\bn C_j^l$ and the chain $\bn D_{>l}+\bn D_0$ cannot contain any cell $\sigma$, such that $i\in A_j^\sigma$ and $\vert A_j^\sigma\vert=l-1$, we conclude that
$\gamma=0$. 

Furthermore, let us fix the sets $A_b$, for all $b\in[p]\setminus \{j\}$.
We denote by $\gamma_{A_1,\dots,\hat{A_j},\dots,A_p}$ the chain of all
cells $\xi=(A_1^\xi,\dots,A_p^\xi)$ from $\gamma$ such that
$A_b^\xi=A_b$, for all $b\in[p]\setminus \{j\}$. Obviously, all these parts
$\gamma_{A_1,\dots,\hat{A_j},\dots,A_p}$ must also be $0$.

Let us also consider the corresponding subchains 
\begin{equation}\label{decomp}
D_{l;\ahat}=\sum_{q\in Q_{l,\ahat}}k_q\sigma_q
\end{equation}
in the chain $D_l$, where 
\[ 
Q_{l,\ahat}=\{q\in Q_l\mid A_b^{\sigma_q}=A_b,\tn{ for all }b\in[p]\setminus\{j\}\}.
\] 
It is clear that $\gamma_{\ahat}$ must be the part of the boundary of $D_{l;\ahat}$, where the deleted element is in $A_j^{\sigma_q}$, and is different from~$i$.

 Let us construct a new cycle  $\widetilde{C}_j^{l}$, homologous to $C_j^l$, with the property that for all cells  $\xi $ appearing in it with $i\in A_j^\xi $ and $A_b^\xi=A_b$, for $b\in[p]\setminus \{j\}$, we have $\vert A_j^\xi  \vert \geq l+1$, that is with the corresponding $D_{l;\ahat}=0$.\\

 $\bullet$ \textsc{Case $l=2$:} Let $\eta$ be a cell from $D_{2;\ahat}$ with coefficient $k\neq 0$, and let $A_j^\eta=\{i,x\}$, for some $x\in[n]\setminus\{i\}$. Since $\gamma_{\ahat}=0$, there must exist another cell $\xi$ in $D_{2;\ahat}$, such that $A_j^\xi=\{i,y\}$, for $y\in[n]\setminus\{i,x\}$. Let us denote $(A_1,\dots,A_{j-1},\{i,x,y\},\dots,A_p)$ with $\sigma$. It is clear that $\sigma\in\Homdm$. Also, the number of cells of the corresponding   $D_{2;\ahat}$ for the cycle $\widehat{C_j^2}=C_j^2-(-1)^K k\; \bn \sigma$, where
\bdm
K=\vert A_1^\eta \vert + \dots+ \vert A_{j-1}^\eta \vert+\left\{ \begin{array}{ll} 1, & \tn{ if } i<y<x\tn{ or }x<y<i; \\ 0, & \tn{otherwise,} \end{array} \right.
\edm  
is reduced at least by one. Repeating this procedure, we will eventually get a chain $\widetilde{C}_j^2$ with the corresponding $D_{2;\ahat}$ equal to zero.\\

$\bullet$  \textsc{Case $l\geq 3$:} Define a map $f_{\ahat}:X_{\ahat}\to \Delta_i$, where 
$X_{\ahat}$ is equal to 
\[\{\eta\in \Homdm\mid A_b^\eta=A_b,\tn{ for }b\in[p] \setminus\{j\},i\in A_j^\eta, \tn{ and }\vert A_j^\eta\vert\geq 2\},
\] 
and $\Delta_i$ is simplex with the vertex set $[n]\setminus\{i\}$, in the following way: for a cell $\eta \in X_{\ahat}$, set
\bdm
f_{\ahat}(\eta)=(-1)^{\vert \{ y \vert y\in A_j^\eta\tn{ and } y>i\}\vert}A_j^\eta \setminus\{i\},
\edm
and then extend it by linearity. Function $f_{\ahat}$ is clearly a bijection between cells and $f_{\ahat}(D_{l;\ahat})$ is a $(l-2)$-chain in $\Delta_i$. Let us now prove that this chain is in fact a cycle. 

For $q\in Q_{l,\ahat}$, let $A_j^{\sigma_q}=[v^q_1,\dots,v^q_{s(q)},i,v^q_{s(q)+1},\dots,v^q_{l-1}]$ where $v^q_1<\dots<v^q_{s(q)}<i<v^q_{s(q)+1}<\dots<v^q_{l-1}$. In order to avoid "ugly" formulas we will use $Q^j_l$ instead of $Q_{l,\ahat}$. Then we have:
\bdm
\begin{split}
\bn (f&_{\ahat}(\sum_{q\in Q^j_l } k_q\sigma_q)) = \sum_{q\in Q^j_l } k_q \bn(f_{\ahat}(\sigma_q))= \\
&=\sum_{q\in Q^j_l } k_q\bn((-1)^{l-1-s(q)}[v^q_1,\dots,v^q_{s(q)},v^q_{s(q)+1},\dots,v^q_{l-1}])= \\ &=
\sum_{q\in Q^j_l } (-1)^{l-1-s(q)} k_q\Big(\sum_{b=1}^{l-1}(-1)^{b-1}[v^q_1,\dots,\hat{v^q_{b}},\dots,v^q_{l-1}]\Big).
\end{split}
\edm
On the other hand we have
\bdm
\begin{split}
0 &=\gamma_{\ahat}= \\ &=(-1)^{\sum_{z=1}^{j-1} \vert A_z \vert}\Big[\sum_{q\in Q^j_l } k_q\Big( \sum_{b=1}^{s(q)}(-1)^{b-1}(A_1^{\sigma_q},\dots,A_j^{\sigma_q}\setminus\{v^q_b\},\dots,A_p^{\sigma_q})+\\ &+ \sum_{b=s(q)+1}^{l-1}(-1)^{b}(A_1^{\sigma_q},\dots,A_j^{\sigma_q}\setminus \{v^q_b\},\dots,A_p^{\sigma_q})\Big)\Big],
\end{split}
\edm
and hence,
\bdm
\begin{split}
0 &=f_{\ahat}(\gamma_{\ahat})=\\ &=(-1)^{\sum_{z=1}^{j-1} \vert A_z \vert} \sum_{q\in Q^j_l } \Big[k_q \Big(\sum_{b=1}^{s(q)}(-1)^{b-1}(-1)^{l-1-s(q)}[v^q_1,\dots,\hat{v^q_{b}},\dots,v^q_{l-1}]+\\ &+ \sum_{b=s(q)+1}^{l-1}(-1)^{b}(-1)^{l-2-s(q)}[v^q_1,\dots,\hat{v^q_{b}},\dots,v^q_{l-1}]\Big)\Big]=\\&=(-1)^{\sum_{z=1}^{j-1} \vert A_z \vert}\sum_{q\in Q^j_l } (-1)^{l-1-s(q)} k_q\Big(\sum_{b=1}^{l-1}(-1)^{b-1}[v^q_1,\dots,\hat{v^q_{b}},\dots,v^q_{l-1}]\Big)=\\&=(-1)^{\sum_{z=1}^{j-1} \vert A_z \vert}\ \bn (f_{\ahat}(\sum_{q\in Q^j_l } k_q\sigma_q)).
\end{split}
\edm
Since a simplex is acyclic, there exists an $(l-1)$-dimensional chain $\tau$ in $\Delta_i$ such that 
\bdm
\bn\tau=f_{\ahat}(\sum_{q\in Q_{l}^j } k_q\sigma_q).
\edm
Let now $\eta_\tau=f^{-1}_{\ahat}(\tau)$. Clearly, $\eta_\tau$ is a~$(t+1)$-chain. We need to check that $\eta_\tau\in \Homdm$. The condition for that is $A_{i_1}^{\eta_\tau}\cap A_{i_2}^{\eta_\tau}=\emptyset$, for $x_{i_1}$ adjacent to $x_{i_2}$. This is clear for $i_1\neq j$, $i_2\neq j$, since then $A_{i_b}^{\eta_\tau}=A_{i_b}$, for $b=1,2$. Assume that there exists $x_{i_1}$ adjacent to $x_j$, such that $A_{i_1}^{\eta_\tau}\cap A_j^{\eta_\tau} \neq\emptyset$. By construction of $\tau$, we know that $A_j^{\eta_\tau}\subseteq\cup_{q\in Q_l^j} A_j^{\sigma_q}$. Since $A_j^{\sigma_q}\cap A_{i_1}=\emptyset$, for any $q\in Q_l^j$, we arrive at a~contradiction. 

Let now $\bn \eta_\tau= \beta_{\tau}+\gamma_{\tau}+\delta_{\tau}$, where 
\begin{eqnarray*}
& &\beta_{\tau} = \{\xi \tn{ from }\bn \eta  _\tau\tn{ such that } 
\vert A_j^\xi  \vert=\vert A_j^{\eta  _\tau} \vert=l+1\}, \\
& &\gamma_{\tau} =\{\xi \tn{ from }\bn \eta  _\tau\tn{ such that } 
\vert A_j^\xi  \vert=l\tn{ and }i\in A_j^\xi \}, \\
& &\delta_{\tau} =\{\xi \tn{ from }\bn \eta  _\tau \tn{ such that } 
i\notin A_j^\xi \}.
\end{eqnarray*}
Similarly to what we have done before, one can prove that 
\bdm
\begin{split} 
f_{\ahat}(\gamma_{\tau}) &= \varepsilon \bn f_{\ahat}(\eta_\tau)=\varepsilon 
\bn \tau= \\ &=f_{\ahat}(\varepsilon D_{l;\ahat}),
\end{split}
\edm
where $\varepsilon=(-1)^{\sum_{z=1}^{j-1}\vert A_z \vert}$. Since $f_{\ahat}$ is a bijection, we conclude that $\gamma_{\tau}=\varepsilon D_{l;\ahat}$. Finally, let $\widetilde{C}_j^{l}=C_j^l-\varepsilon\bn \eta  _\tau$.\\

After repeating the procedure described above for all combinations of sets $A_1,\dots,A_{j-1},A_{j+1},\dots,A_p$ appearing in the decomposition of $D_l$ to subchains of the form (\ref{decomp}), we will get a chain $C_j^{l+1}$ such that  $\vert A_j^\xi  \vert \geq l+1$, for all cells  $\xi $ appearing in it with nonzero coefficients, and with $i\in A_j^\xi $. 
 
The following observation is important for our argument. If for all cells $\xi $ from $C_j$ we have that $A_b^\xi  \subseteq S_b$, for $b\in [p]$ and $S_b\subseteq [n]$, then also for all $\xi'$ from $C_{j+1}$ and $b\in [p]$, we have $A_b^{\xi'} \subseteq S_b$.

 After first phase, we shall eventually obtain a~chain $C_1:=C_{p+1}$ which has the following additional property: for all cells $\eta$ appearing in this chain with a~nonzero coefficient, we have $A_j^\eta \subseteq\{i-1,i,\dots,n\}$, for $j\in[\lambda]$ (the original condition of the lemma preserved), and $i\notin A_j^\eta$, for $j\in \{\lambda+1,\dots,p\}$. %Denote $C_{p+1}$ with $C_1$ and repeat the following two steps for all $j\in [\lambda]$

\nin {\bf Step 3.} This step is almost the same as Step 1. Namely, assume $\eta$ is a~cell from $C_j$ such that coefficient of $\eta$ in $C_j$ is $k$, and $A_j^\eta=\{i-1\}$. Then $\eta'=(A_1^\eta,\dots,A_{j-1}^\eta,\{i-1,i\},A_{j+1}^\eta,\dots,A_p^\eta)\in \Homdm$, since $i\notin A_q^\eta$, for all $q$ such that $x_q\in \texttt{N}(x_j)$. Again, $C'_{j}=C_j-(-1)^{K(\eta')}k \ \bn \eta'$, where $K(\eta')=\vert A_1^\eta \vert + \dots+ \vert A_{j-1}^\eta \vert+1$, 
is a chain homologous to $C_j$ and number of cells $\sigma$ from $C_j'$ such that $A_j^\sigma=\{i-1\}$ is strictly less then number of such cells from $C_j$. We repeat this procedure until we get cycle $C_j^2$ in which there are no cells  $\sigma$ such that $A_j^\sigma=\{i-1\}$.

\nin {\bf Step 4.} This step is very similar to the Step 2. The only difference is that we are "removing" $i-1$ (instead of $i$) from $j$-th coordinate sets of each cell from the chain $C_j$ and obtain a new chain $C_{j+1}$ such that $i-1 \notin A_j^\eta$ for each $\eta$ from $C_{j+1}$. By the observation that we have made after Step 2, we see that, for all cells $\eta$ from $C_{j+1}$ and for all $q\in [j]$, we have $A_j^\eta \subseteq\{i,i+1,\dots,n\}$. In the case when $j=\lambda$, we set $C':=C_{j+1}$.

Now it is easy to see that $C'\in X_i(\lambda)$ and that $C-C'$, is a~boundary of some $(t+1)$-chain of \Homk, and we have proved our claim.
\end{proof}
\end{lem}

\begin{figure}[ht] 
\begin{center}
\psfrag{1}{$\scriptstyle{1}$} 
\psfrag{2}{$\scriptstyle{2}$} 
\psfrag{3}{$\scriptstyle{3}$} 
\psfrag{4}{$\scriptstyle{4}$} 
\psfrag{a}{$\scriptstyle{(1,6,235,2)}$}
\psfrag{b}{$\scriptstyle{(1,67,25,2)}$}
\psfrag{c}{$\scriptstyle{(1,67,35,2)}$}
\psfrag{5}{$\scriptstyle{(1,7,235,2)}$}
\psfrag{6}{$\scriptstyle{(1,7,234,2)}$}
\psfrag{7}{$\scriptstyle{(1,67,34,2)}$}
\psfrag{8}{$\scriptstyle{(1,67,24,2)}$}
\psfrag{9}{$\scriptstyle{(1,6,234,2)}$}
\psfrag{(a)}{(a)}
\psfrag{(b)}{(b)}
%\psfrag{1}{$\scriptscriptstyle{(1,2,1,2)}$}
%\psfrag{Hom2}{$\Hom (C_4,C_{2r})$}
\includegraphics[scale=0.9]{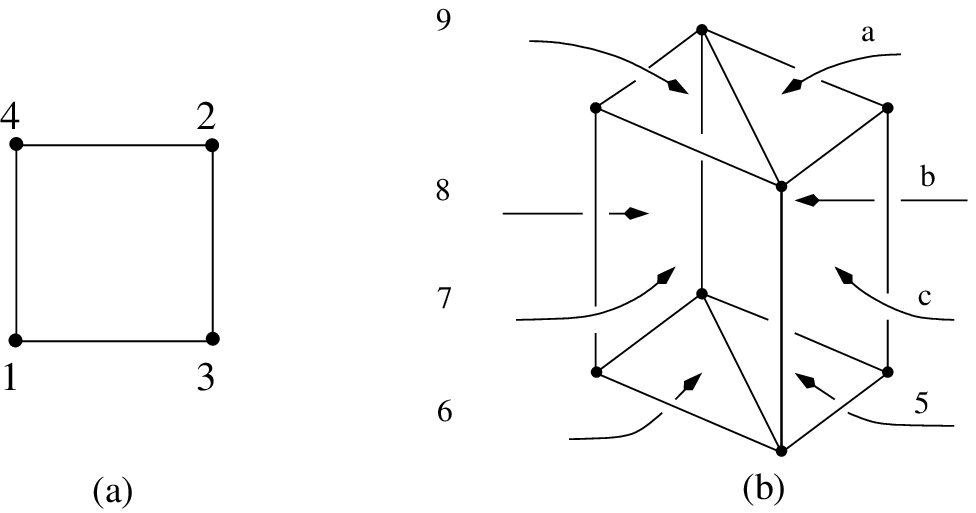}
\caption{} 
\label{example}
\end{center} 
\end{figure}

\textbf{Example.} Now we will illustrate the proof of Lemma \ref{dokaz_hom} on a cycle from $\texttt{Hom}\,(C_4,K_7)$. For simplicity, we will use slightly different notations for cells, for example we will write $(1,67,234,2)$ instead of $(\{1\},\{6,7\},\{2,3,4\},\{2\})$. We have ordered vertices of $C_4$ as described in  Lemma \ref{dokaz}, see Figure \ref{example}(a). 

Let $C=C_3=-(1,7,234,2)+(1,6,234,2)-(1,67,34,2)+(1,67,24,2)+(1,7,235,2)-(1,6,235,2)+(1,67,35,2)-(1,67,25,2)\in X_1(2)$. It is easy to check that $\bn C=0$. Figure \ref{example}(b) depicts the part of the complex  $\texttt{Hom}\,(C_4,K_7)$ formed by all the
cells from $C$.

 First we will eliminate $2$ from $A_3^\eta$ and $A_4^\eta$, for all cells $\eta$ from our chain $C$. 

For $j=3$ we obviously do not have Step 1. Hence, $C=C_3^2=D_{2;1,67,\hat{A}_3,2}+D_{>2}+D_0$, where $D_{2;1,67,\hat{A}_3,2}=(1,67,24,2)-(1,67,25,2)$, $D_{0}=(1,67,35,2)-(1,67,34,2)$, and $D_{>2}$ is the rest of $C_3^2$. Let $\eta=(1,67,24,2)$ and $\xi=(1,67,25,2)$. Then $\sigma=(1,67,245,2)$, $k=+1$, $K=3$ and finally, 
\bdm
\begin{split}
C_3^3 =&\widetilde{C}_3^2=C_3^2+\bn \sigma =(1,6,234,2)-(1,7,234,2)-(1,67,34,2)+(1,7,235,2)-\\ &-(1,6,235,2)+(1,67,35,2)-(1,7,245,2)+(1,6,245,2)-(1,67,45,2).
\end{split}
\edm
Again, $D_{3;1,7,\hat{A}_j,2}=-(1,7,234,2)+(1,7,235,2)-(1,7,245,2)$, $D_{3;1,6,\hat{A}_j,2}=(1,6,234,2)-(1,6,235,2)+(1,6,245,2)$ and $D_0=C_3^3-D_{3;1,7,\hat{A}_j,2}-D_{3;1,6,\hat{A}_j,2}$, since $D_{>3}=0$. Then $f_ {3;1,7,\hat{A}_j,2}(D_{3;1,7,\hat{A}_j,2})=-(34)+(35)-(45)=\bn \tau$, where $\tau=-(345)$. It is easy to see that $\eta_\tau=(1,7,2345,2)$ and $\widetilde{C}_3^3=C_3^3-\bn\eta_\tau=(1,6,234,2)-(1,67,34,2)-(1,6,235,2)+(1,67,35,2)+(1,6,245,2)-(1,67,45,2)-(1,7,345,2)$.  If we do the same thing for the subchain $D_{3;1,6,\hat{A}_j,2}$, we will get the cycle $C_4=-(1,67,34,2)+(1,67,35,2)-(1,67,45,2)-(1,7,345,2)+(1,6,345,2)$. 

Now we apply Step 1 for $j=4$: Let $\eta=(1,67,34,2)$, $k=-1$. We see that, in this case, $I(\eta)$ could be any element from the set $\{3,4,5\}$, but for convenience, we will choose $\eta'=(1,67,34,23)$. Then $C_4'=C_4+\bn\eta'$.  Repeating this for all cells such that their fourth coordinate set is $\{2\}$, with additional remark that we always pick $I(\eta)=3$, we get a chain $C_1=-(1,67,34,3)+(1,67,35,3)-(1,7,345,3)+(1,6,345,3)-(1,67,45,3)$. 

Let us now remove $1$ from the first two coordinate sets of all cells from $C_1$.

 First we apply Step 3: Let $\eta=(1,67,34,3)$. Then $\eta'=(12,67,34,3)$ and $C_1'=C_1-\bn \eta'$. If we do the same thing for all other cells from $C_1$, finally we get $C'=-(2,67,34,3)+(2,67,35,3)-(2,67,45,3)-(2,7,345,3)+(2,6,345,3)\in X_2(2)$.
\begin{thm}\label{homology} 
Let $G$ be an arbitrary graph, then
\textnormal{\bdm H_t(\Homdm,\mathbb{Z})=0,\textnormal{\it{ for all }} 1\leq t\leq \vgap(G,n)-1.
\edm}
\nin{\bf Note.}
As mentioned above, it was proved in \cite{BK03b} that \emph{$H_0(\Homdm,\mathbb{Z})=\mathbb{Z}$} if $\vgap(G,n)\geq 1$. Also, it is a~direct corollary of Theorem \ref{fgrupa} that \emph{$H_1(\Homdm,\mathbb{Z})=0$} if $\vgap(G,n)\geq 2$. Hence, we can assume that $t\geq 2$.

\begin{proof}
Let $p$ and $\lambda$ be the same as in Lemma \ref{dokaz} and Lemma \ref{dokaz_hom}. Also, let $C$ be a $t$-cycle from \Homk, that is dimension of all cells from $C$ is equal to $t$ and $\bn C=0$, where the boundary operator is defined in the proof of Lemma \ref{dokaz_hom}. In order to prove that $H_t(\Homdm,\mathbb{Z})=0$, we need to prove that $C$ bounds. Again, we will use induction on the maximal valency of $G$.

Let the maximal valency of a graph $G$ be zero. In the Theorem \ref{fgrupa} we proved that in this case \Hom is contractible and hence, $H_t(\Homdm,\mathbb{Z})=0$ for all $0\leq t \leq n-2$ and there exists a chain $D$ such that $C=\bn D$. 

Assume now that the maximal valency of $G$ is $d\geq 1$. By Lemma \ref{dokaz_hom}, there exists a cycle $C'$ which is homologous to $C$ and which has the property that, for all cells $\eta$ with nonzero coefficients in this cycle and for all $q\in [\lambda]$, $A_q^\eta=\{n\}$. Hence, $C'$ is isomorphic to a cycle inside $\texttt{Hom}\,(G-\{x_1,x_2,\dots,x_\lambda\},K_{n-1})$. The set $\{x_1,x_2,\dots,x_\lambda\}$ is a maximal independent set in $G$, therefore, the maximal valency of $G-\{x_1,x_2,\dots,x_\lambda\}$ is strictly less then the maximal valency of $G$. It follows that $C'$ is a~boundary in $\texttt{Hom}\,(G-\{x_1,x_2,\dots,x_\lambda\},K_{n-1})$ by the induction hypothesis.  We conclude that there exists a $(t+1)$-chain $D$ in $\thom(G,K_n)$, such that $\bn D=C$.
\end{proof}
\end{thm}

\begin{rem} It was proved in \cite{BK03a} that:
\begin{center}$H_i(\texttt{Hom}\, (C_{2r+1},K_n);\mathbb{Z}) = 0$, for $1\le i\le n-4$.
\end{center} 
 This fact is a direct corollary of the previous theorem since the maximal valency of $C_{2r+1}$ is equal to $2$. 
\end{rem}

Finally, we are able to put the pieces together and prove the Conjecture~\ref{conj:main}.

%\begin{thm}\label{THE} 
%The following inequality is valid for an arbitrary graph $G$:
%\begin{equation} \label{eq:main}
 %\conn\thom(G,K_n)\geq\vgap(G,n)-1.
%\end{equation}

\pet

\nin {\em Proof of the Conjecture~\ref{conj:main}.} 

\pet

\nin The cases when $\vgap(G,n)=0,1$ were observed in \cite[Proposition
2.4]{BK03b}. The case $\vgap(G,n)=2$ is proven in Theorem
\ref{fgrupa}.

Let us now deal with case $\vgap(G,n)\geq 3$. By the Theorem \ref{fgrupa}, \Hom is 1-connected, and by the Theorem \ref{homology} we have,
\bdm
 H_t(\Homdm,\mathbb{Z})=H_t(\Homdm)=0,
\edm
 for all  $1\leq t\leq\vgap(G,n)-1$. By a~standard corollary to the Hurewitz theorem (see, for example, \cite[10.10 Corollary, page 479]{Bre}), we have that $\pi_t(\Homdm,*)=0$, for all $t\leq\vgap(G,n)-1$, and, hence $\conn\thom(G,K_n)\geq\vgap(G,n)-1$ by definition.
\qed
\begin{rem} We know that the result of Conjecture~\ref{conj:main} is sharp for several classes of graphs, for example for odd cycles and complete graphs.
\end{rem}

\begin{rem} By the same corollary we used in the proof of Conjecture \ref{conj:main}, we conclude that $\pi_{\vgap(G,n)}(\Homdm,*)\approx H_{\vgap(G,n)}(\Homdm)$.
\end{rem} 

\begin{cor} The complex \emph{$\texttt{Hom}\,(C_r,K_{n})$} is $(n-4)$-connected, for arbitrary integers $r,n\geq 3$.
\end{cor}

\end{document}